\documentclass{amsart}

\newtheorem{theorem}{Theorem}
\newcommand{\bt}{\begin{theorem}}
\newcommand{\et}{\end{theorem}}

\newtheorem{lemma}{Lemma}
\newcommand{\bl}{\begin{lemma}}
\newcommand{\el}{\end{lemma}}

\newcommand{\N}{\ensuremath{ \mathbf{N} }}
\newcommand{\Z}{\ensuremath{ \mathbf{Z} }}

\DeclareMathOperator{\qqand}{\qquad\text{and}\qquad}

\usepackage{amssymb,latexsym}

\title{Congruence classes and maximal nonbases}
\author{Melvyn B. Nathanson}
\address{Department of Mathematics\\Lehman College (CUNY)\\Bronx, NY 10468}
 \email{melvyn.nathanson@lehman.cuny.edu}
 \thanks{This research was supported by a grant from the PSC-CUNY Research Awards Program.}
 \dedicatory{To Endre Szemer{\' e}di on his 80th birthday}



\date{\today}
 
\begin{document}
 
 \begin{abstract}
The set $A$ is an asymptotic nonbasis of order $h$ for an additive abelian semigroup $X$ 
if there are infinitely many elements of $X$ not in the $h$-fold sumset $hA$.
For all $h \geq 2$, this paper constructs new classes of asymptotic nonbases of order $h$ 
for \Z\ and for $\N_0$  
that are not subsets of maximal asymptotic nonbases.
\end{abstract}

 \maketitle
 
\section{Asymptotic bases and nonbases}

Let $X$ be an additive abelian group or semigroup,  
and let $A$ be a subset of $X$.  
For every positive integer $h$, we define the \emph{$h$-fold sumset }
\[
hA =  \{a_1 + \cdots + a_h: a_i \in A \text{ for all } i = 1,\ldots, h \}.
\]
The subset $A$  is an \emph{basis of order $h$ for $X$} if every element of $X$ can be represented 
as the sum of $h$ not necessarily distinct elements of $A$, 
that is, if $hA = X$.  
The subset $A$ is a \emph{nonbasis of order $h$ for $X$} if $hA \neq X$.  
The nonbasis $A$ is a \emph{maximal nonbasis of order $h$  for $X$}\index{maximal nonbasis} 
if $A\cup \{b\}$ is an basis of order $h$ for $X$ for all $b \in X\setminus A$.

The subset $A$ of a set $X$ is 
\emph{co-finite}\index{co-finite} if $X\setminus A$  is finite and 
\emph{co-infinite}\index{co-infinite} if $X\setminus A$  is infinite.  

The subset $A$ of an infinite semigroup $X$  is an \emph{asymptotic basis of order $h$ for $X$} 
if the sumset $hA$ contains all but finitely many elements of $X$, 
that is, if $hA$ is co-finite, and an 
\emph{asymptotic nonbasis of order $h$}\index{asymptotic nonbasis} 
if $A$ is not an asymptotic basis of order $h$ for $X$, that is, if $hA$ is co-infinite.
An asymptotic nonbasis $A$  of order $h$  for $X$ is a 
\emph{maximal}\index{maximal asymptotic nonbasis} 
if $A\cup \{b\}$ is an asymptotic basis of order $h$ for $X$ 
for all $b \in X\setminus A$.

The subset $A$ of an infinite semigroup $X$ is an asymptotic nonbasis of order 1 
for $X$ if and only if it is co-infinite.  
Because there is no maximal co-infinite subset of an infinite set, 
there exists no maximal asymptotic nonbasis of order 1 for $X$.   
In particular, no asymptotic nonbasis of order 1 for $X$ is a subset of 
a maximal asymptotic nonbasis of order 1 for $X$.

In this paper we consider the additive group \Z\ of integers 
and the additive semigroup $\N_0$ of nonnegative integers.   
There is no maximal asymptotic nonbasis of order 1 for \Z\ or for $\N_0$.
For $h \geq 2$, there do exist maximal asymptotic nonbases of order $h$ for \Z\ and for $\N_0$.
The first examples were constructed in~\cite{nath1974-11}.  
It had been an open problem 
to determine if every asymptotic nonbasis of order $h \geq 2$ is a subset of a 
maximal asymptotic nonbasis of order $h$.
Hennefeld~\cite{henn77} constructed the first example of an asymptotic nonbasis of order $h$ 
that could not be embedded in a maximal asymptotic nonbasis of order $h$.
Recently, Ling~\cite{ling17} obtained two new classes of asymptotic nonbases of order $h$  
that are not contained in maximal asymptotic nonbases of order $h$.
This paper describes a simple construction that includes 
Hennefeld's and Ling's examples as special cases.

Nonbases and maximal nonbases were introduced by Nathanson~\cite{nath1974-11}, 
and investigated 
by Deshouillers and Grekos~\cite{desh-grek79}, 
Erd\H os and Nathanson~\cite{nath1975-13, nath1975-17,nath1976-22, nath1977-29,nath1979-36}, 
Hennefeld~\cite{henn77}, Ling~\cite{ling17,ling18}, 
Nathanson~\cite{nath1977-28}, and 
Nathanson and S\' ark\" ozy~\cite{nath1996-87}.  For surveys of open problems, see 
Erd\H os and Nathanson~\cite{nath1987-59} 
and Nathanson~\cite{nath1989-71,nath2009-134}.

\section{Asymptotic nonbases for \Z}
For $d \in \Z$ and $A \subseteq \Z$, define the \emph{dilation}\index{dilation} 
$d\ast A = \{da:a\in A\}$.

\bt                  \label{MaximalNonbases:theorem:CongruenceBasis}
Let $h\geq 2$.  Let $s$ and $t$ be integers, and let  
\[
A_{\Z} = \{s\} \cup \{ hz +t: z \in \Z \}.
\]
If 
\begin{equation}                  \label{MaximalNonbases:gcd-d}
\gcd(h,s-t) = d \geq 2
\end{equation}
then $hA_{\Z} \subseteq d\ast \Z$ and $A_{\Z}$ is an asymptotic  nonbasis of order $h$ for \Z. 

If 
\begin{equation}                  \label{MaximalNonbases:gcd}
\gcd(h,s-t) = 1
\end{equation}
then $A_{\Z}$ is a  basis of order $h$ for \Z. 
Moreover, if $n  \equiv t-s  \pmod{h}$, then there exists $z_1\in \Z$ such that 
\begin{equation}                  \label{MaximalNonbases:s-t}
n = (h-1)s + \left( hz_1  + t \right) \in hA_{\Z} 
\end{equation}
and this is the unique representation of $n$ as a sum of $h$ elements of $A_{\Z}$.  
\et
 
\begin{proof}
Let $n \in hA_{\Z}$.  There exist $i \in \{0,1,\ldots, h\}$ and $z_1,\ldots, z_{h-i} \in \Z$ 
such that 
\[
n = is + \sum_{j=1}^{h-i} (hz_j+t) \equiv is+ (h-i)t \equiv i(s-t) \pmod{h}.
\]
If $h$, $s$, and $t$ satisfy the divisibility condition~\eqref{MaximalNonbases:gcd-d}, then 
$n \equiv 0 \pmod{d}$, and so $hA_{\Z} \subseteq d\ast \Z$.  
It follows that $A_{\Z}$ is an asymptotic nonbasis of order $h$ for \Z.

If $h$, $s$, and $t$ satisfy the divisibility condition~\eqref{MaximalNonbases:gcd},  
then 
\[
i(s-t) \not\equiv j(s-t)\pmod{h} 
\]
for all $i,j \in \{0,1,\ldots, h-1\}$ with $i \neq j$.  
It follows that, for all $n \in \Z$, there are unique integers 
\[
i \in \{0,1,\ldots, h-1\} \qqand q \in \Z
\]
such that 
\[
n \equiv i(s-t) \pmod{h} 
\]
and  
\[
n = i(s-t) + hq. 
\]
Let 
\[
k = h-i \in \{1,2,\ldots, h\}
\]
and let $z_1,\ldots, z_k$ be integers such that 
\[
\sum_{j=1}^k z_j = q - t.
\]
We have 
\begin{align}
n & =  i(s-t) + hq       \nonumber  \\
& = is + h\left( q - t \right) + kt       \nonumber  \\
& = is + h\left( \sum_{j=1}^k z_j \right) + kt      \nonumber  \\
& = is +  \sum_{j=1}^k \left( h z_j  + t \right) \in hA      \label{MaximalNonbases:nis}
\end{align}
and so $A_{\Z} $ is a basis of order $h$ for \Z.

For all $n \in \Z$, if $i \in \{0,1,\ldots, h-1\}$,  $k = h-i$,  
and $z_1,\ldots, z_k \in \Z$ satisfy  
\[
n =  is +  \sum_{j=1}^k \left( h z_k + t \right) 
\]
then 
\[
n \equiv is + kt \equiv i (s-t)  \pmod{h}
\]
and so $n$ uniquely determines the integer $i$.  
Moreover, $i = h-1$ if and only if $n  \equiv t-s \pmod{h}$, 
and~\eqref{MaximalNonbases:s-t} is the unique representation of $n$ as the sum of $h$ elements of $A_{\Z}$.  
This completes the proof.  
\end{proof}

\bl                           \label{MaximalNonbases:lemma:Ygap}  
Let $Y$ be a set of integers such that 
there exist only finitely many pairs $y,y' \in Y$ 
with $y \neq y'$ and  $|y-y'| \leq 3$.
The set $X = \Z\setminus Y$ is an asymptotic basis of order $h$ for \Z\ 
for all $h \geq 2$.

Let $Y_0$ be a set of  nonnegative integers such that 
there exist only finitely many pairs $y,y' \in Y_0$ 
with $y \neq y'$ and  $|y-y'| \leq 3$.
The set $X_0 = \N_0 \setminus Y_0$ is an asymptotic basis of order $h$ for $\N_0$
for all $h \geq 2$.
\el

\begin{proof}
Let $n \in \Z$.  If $n  = 2u$ is even, then 
\[
n =  u + u= (u+1) + (u -1).
\]
 If $n = 2u+1$ is odd, then 
\[
n = u + (u +1) =  (u-1) + (u + 2).
\]
The gap condition $|y-y'| \leq 3$ implies that 
\begin{equation}                                       \label{MaximalNonbases:Ygap}
|\{u-1,u,u+1,u+2\} \cap Y| \geq 2
\end{equation}
for only finitely many integers $u$, 
and so $n \in 2X$ for all but at most finitely many integers $n$.
Thus, $X$ is an asymptotic basis of order 2 for \Z, 
and there is a finite subset $F$ of \Z\ such that $2X = \Z\setminus F$.
  
Let $h \geq 3$.  Choose  $x_0 \in X$.  
If $n -(h-2)x_0 \notin F$, then there exist $x,x' \in X$ such that 
\[
n -(h-2)x_0 = x + x'
\]
and so $n = (h-2)x_0 + x + x' \in hX$.
Therefore, $X$ is an asymptotic basis of order $h$ for \Z\ for all $h \geq 2$.

The proof for $\N_0$ is similar.  
\end{proof}

The set $Y$  of integers has \emph{infinite gaps}\index{infinite gaps} 
if, for all $C > 0$, 
there exist only finitely many pairs of integers $y,y' \in Y$ 
with $y \neq y'$ and $|y-y'| \leq C$.

\bt                        \label{MaximalNonbases:theorem:CongruenceNonbasis}
Let $h \geq 2$.    
Let $s$ and $t$ be integers such that 
\begin{equation}                                                  \label{MaximalNonbases:CongruenceNonbasis-11}
\gcd(h,s-t) = 1.
\end{equation}
Let $Y$ be an infinite set of integers with infinite gaps, and let
\[
X = \Z \setminus Y. 
\]
The set 
\[
A_X = \{s\} \cup \{hx+t: x \in X \}
\]
is an asymptotic nonbasis of order $h$ for \Z\ 
that is not a subset of a maximal asymptotic nonbasis of order $h$ for \Z. 
\et

\begin{proof}
We begin by proving that there is a finite set $\mathcal{F}$ of integers such that 
\begin{equation}                                                   \label{MaximalNonbases:CongruenceNonbasis-hA}
\Z \setminus hA_X = \mathcal{F} \cup \{ (h-1)s + hy+t : y \in Y \}
\end{equation}
and
\begin{equation}                                                   \label{MaximalNonbases:CongruenceNonbasis-hA-F}
\mathcal{F} \cap \{ (h-1)s + hy+t : y \in Y \} = \emptyset.
\end{equation}

Let $n \in \Z$ satisfy 
\begin{equation}                                                  \label{MaximalNonbases:CongruenceNonbasis-22}
n \not\equiv t-s \pmod{h}.
\end{equation}
The divisibility condition~\eqref{MaximalNonbases:CongruenceNonbasis-11} 
and the congruence condition~\eqref{MaximalNonbases:CongruenceNonbasis-22} imply that 
there is a unique integer $i \in \{0,1,\ldots, h-2\}$ such that 
\[
n \equiv i(s-t)\pmod{h}.
\]
Let $k = h-i \in \{2,3 \ldots, h\}$.  There exists a unique integer $q$ such that 
\[
n = i(s-t) + hq = is + h(q - t) +kt.
\]
By Lemma~\ref{MaximalNonbases:lemma:Ygap}, 
for all but finitely many integers $q$, there exist $x_1,\ldots, x_k \in X$ such that 
\[
x_1+\cdots + x_k = q-t
\]
and so 
\begin{align*}
n & = is + h(q - t) +kt \\
& = is + h \left(x_1+\cdots + x_k \right)  +kt \\
& = is + \sum_{j=1}^k (hx_j +t) \in hA_X 
\end{align*} 
for all but finitely many integers $n \not\equiv t-s \pmod{h}$, 

Let $\mathcal{F}$ be the finite set of integers $n$ such that $n \not\equiv t-s \pmod{h}$ 
and $n \notin hA_X$.

Let $A_{\Z} = \{s\} \cup \{hz+t: z\in \Z\}$.  
By Theorem~\ref{MaximalNonbases:theorem:CongruenceBasis}, 
if $n \in \Z$ and 
\[
n \equiv t-s \pmod{h}
\]
then there is a unique integer $z_1 \in \Z$ such that 
\[
n = (h-1)s + (hz_1+t)
\]
and this is the unique representation of $n$ as a sum of $h$ elements of $A_{\Z}$.
If $z_1 = x_1\in X$, then $hx_1+t \in A_X$ and $n \in hA_X$. 
If $z_1 = y_1 \in Y$, then $hy_1+t \notin A_X$ and $n \notin hA_X$.  
This proves~\eqref{MaximalNonbases:CongruenceNonbasis-hA}.  
The set $Y$ is infinite, 
and so $A_X$ is an asymptotic nonbasis of order $h$ for \Z.  

We shall prove that $A_X$ is not a subset of a maximal asymptotic nonbasis of order $h$ for \Z.   

Let  $b \in  \Z\setminus A_X$ and  
\[
b \not\equiv s,t \pmod{h}.
\]
Because $\gcd(h, s-t) = 1$, the congruence 
\begin{equation}                                                \label{MaximalNonbases:CongruenceNonbasis-ist}
(i+1)(s-t) \equiv t-b \pmod{h} 
\end{equation}                                            
has a unique solution $i\in \{0,1,2,\ldots, h-1\}$.
If $i = h-1$, then $b \equiv t \pmod{h}$, which is absurd.
If $i = h-2$,  then $b \equiv s \pmod{h}$, which is also absurd. 
Therefore,  $i \in \{0,1,2,\ldots, h-3 \}$ and $h-i-1 \in \{2,3,\ldots, h-1\}$.  

By~\eqref{MaximalNonbases:CongruenceNonbasis-hA} 
and~\eqref{MaximalNonbases:CongruenceNonbasis-hA-F}, 
we have $n \notin hA_X \cup \mathcal{F}$ 
if and only if there exists $y \in Y$ such that 
\[
n = (h-1)s + hy + t.
\]
Let $i \in \{0,1,2,\ldots, h-3 \}$ satisfy~\eqref{MaximalNonbases:CongruenceNonbasis-ist}.  
There exists a unique integer  $w$ such that 
\begin{equation}                                                \label{MaximalNonbases:CongruenceNonbasis-case2}
n = (h-1)s + hy + t =  b + is + hw  + (h-i-1)t. 
\end{equation}
Moreover,  $h-i-1\geq 2$.    
By Lemma~\ref{MaximalNonbases:lemma:Ygap}, the set $X$ is an asymptotic basis of order $h-i-1$.  
For all but finitely many $w$,  
there exist integers $x_1,x_2,\ldots, x_{h-i-1} \in X$ such that 
\[
w = \sum_{j=1}^{h-i-1} x_j
\]
and so 
\begin{align*}
n & =  b+ is  +  hw + (h-i-1)t \\
& =  b+ is  + h \sum_{j=1}^{h-i-1} x_j  + (h-i-1)t \\
& =  b+ is  +  \sum_{j=1}^{h-i-1} ( hx_j +t) \\
& \in h \left( A_X \cup \{b\} \right).
\end{align*}
Therefore, if $b \in  \Z\setminus A_X$ and $b \not\equiv s,t \pmod{h}$, then $A_X \cup \{b\} $ 
is an asymptotic basis of order $h$.

Let $b \in  \Z\setminus A_X$ and 
\[
b \equiv s \pmod{h}.
\]
We have  $b \neq s$ and $b = s+hu$ for some $u \neq 0$.
For all $y \in Y$, 
\begin{align*}
(h-1)s+ hy + t  & = (h-1)(s+hu) + h(y-(h-1)u) + t \\ 
 & = (h-1)b + h(y-(h-1)u) + t.
\end{align*}
Because $Y$ has infinite gaps, there exist only finitely many integers $y \in Y$ 
with $y-(h-1)u \in Y$, 
and so $y-(h-1)u \in X$ for all but finitely many $y \in Y$. 
It follows that $(h-1)s+ hy + t \in h( A_X \cup \{b\}  )$ for all but finitely many $y \in Y$. 
Therefore, if $b \in  \Z\setminus A_X$ and $b \equiv s \pmod{h}$, then $A_X \cup \{b\} $ 
is an asymptotic basis of order $h$.

Let $b \in  \Z\setminus A_X$ and 
\[
b \equiv t \pmod{h}.
\]
We have  
\[
b = hy' + t \qquad \text{for some $y' \in Y$.} 
\]
From~\eqref{MaximalNonbases:CongruenceNonbasis-hA} 
and~\eqref{MaximalNonbases:CongruenceNonbasis-hA-F}, we have 
\[
\Z \setminus \left( hA_X \cup  \mathcal{F} \right) = \{ (h-1)s + hy+t : y \in Y \}.
\]
The uniqueness statement~\eqref{MaximalNonbases:s-t} in  Theorem~\ref{MaximalNonbases:theorem:CongruenceBasis}  
implies that 
\[
h \left( A_X \cup \{b\} \right) =  hA_X \cup  \mathcal{F'} \cup \{(h-1)s+b\} 
\]
for some $\mathcal{F'} \subseteq \mathcal{F}$.  
Therefore, if $B \subseteq \Z\setminus A_X$ and $A_X \cup B$ is an asymptotic nonbasis 
of order $h$, then there exists $Y' \subseteq Y$ such that 
\[
B = B_{Y'} = \{hy' +t:y' \in Y'\} 
\]
and 
\[
h(A_X\cup B_{Y'}  ) = hA_X \cup  \mathcal{F''} \cup \{ (h-1)s+hy' +t:y' \in Y'\}   
\] 
for some $\mathcal{F'} \subseteq \mathcal{F}$.  
We conclude that $A_X\cup B_{Y'}  $ is an asymptotic nonbasis of order $h$ 
for \Z\ if and only if $Y'$ is a co-infinite subset of $Y$.
Because the infinite set $Y$ contains no maximal co-infinite subset, 
there exists no set $Y'$ such that $A_X \cup B_{Y'} $ is a maximal asymptotic nonbasis of order $h$. 
Thus, the set $A_X$ is an asymptotic nonbasis of order $h$ 
for \Z\ that is not contained in a maximal asymptotic nonbasis of order $h$ for \Z.
This completes the proof.  
\end{proof}

\section{Asymptotic nonbases for $\N_0$}
The proofs of Theorems~\ref{MaximalNonbases:theorem:CongruenceBasis-N0} 
and~\ref{MaximalNonbases:theorem:CongruenceNonbasis-N0} are similar to 
those of Theorems~\ref{MaximalNonbases:theorem:CongruenceBasis}
and~\ref{MaximalNonbases:theorem:CongruenceNonbasis}, 
but with a few subtle differences.

\bt                  \label{MaximalNonbases:theorem:CongruenceBasis-N0}
Let $h\geq 2$.  Let $s$ and $t$ be nonnegative integers, and let  
\[
A_{\N_0} =  \{s\} \cup \{ hz +t: z \in \N_0 \}
\]
If 
\begin{equation}                  \label{MaximalNonbases:gcd-d0}
\gcd(h,s-t) = d \geq 2
\end{equation}
then $hA_{\N_0} \subseteq d\ast \N_0$ and $A$ is an asymptotic  nonbasis of order $h$ for $\N_0$. 

If 
\begin{equation}                  \label{MaximalNonbases:gcd01}
\gcd(h,s-t) = 1
\end{equation}
then $A_{\N_0}$ is an asymptotic  basis of order $h$ for $\N_0$. 
Moreover, if  $n \in \N_0$ and $n  \equiv t-s  \pmod{h}$, then there exists $z_1\in \N_0$ such that 
\begin{equation}                  \label{MaximalNonbases:s-t-N0}
n = (h-1)s + \left( z_1 h + t \right) \in hA_{\N_0} 
\end{equation}
and this is the unique representation of $n$ as a sum of $h$ elements of $A_{\N_0}$.  
\et

\begin{proof}
The proof that $\gcd(h,s-t)  \geq 2$ implies that  $A$ is an asymptotic  nonbasis of order $h$ for $\N_0$ 
is the same as in Theorem~\ref{MaximalNonbases:theorem:CongruenceBasis}.

Suppose that $h$, $s$, and $t$ satisfy the divisibility condition~\eqref{MaximalNonbases:gcd01}.  
For all $n\in \Z$, there are unique integers $i \in \{0,1,\ldots, h-1\}$
and $q$ such that 
\[
n =i(s-t)+ hq. 
\]
If 
\[
n \geq (h-1) |s-t| + ht 
\]
then 
\[
q - t = \frac{n -i(s-t) - ht}{h} \geq \frac{n - (h-1) |s-t| - ht}{h} \geq 0.
\]
Let 
\[
k = h-i \in \{1,2,\ldots, h\}.  
\]
Because $q-t \geq 0$, 
there exist nonnegative integers $z_1,\ldots, z_k$ such that 
\[
\sum_{j=1}^k z_j =  q - t.
\]
We have 
\begin{align*}
n & =  i(s-t)  + hq       \\
& = is + h\left(  q - t \right) + kt       \\
& = is + h\left( \sum_{j=1}^k z_j \right) + kt        \\
& = is +  \sum_{j=1}^k \left( h z_j  + t \right) \in hA     
\end{align*}
and so $A$ is an asymptotic basis of order $h$ for $\N_0$.

Moreover, for all $n \geq (h-1) |s-t|$, if $i \in \{0,1,\ldots, h-1\}$,  $k = h-i$,  
and $z_1,\ldots, z_k \in \N_0$ satisfy  
\[
n =  is +  \sum_{j=1}^k \left( h z_k + t \right) 
\]
then 
\[
n \equiv is + kt \equiv i (s-t)  \pmod{h}
\]
and so $n$ uniquely determines the integer $i$.  
In particular, if $i = h-1$ and $k = 1$, then $n  \equiv t-s \pmod{h}$
and~\eqref{MaximalNonbases:s-t-N0}  
is the unique representation of $n$ as the sum of $h$ elements of $A_{\N_0}$.  
This completes the proof.  
\end{proof}

\bt                        \label{MaximalNonbases:theorem:CongruenceNonbasis-N0}
Let $h \geq 2$.    
Let $s$ and $t$ be nonnegative integers such that 
\begin{equation}                               \label{MaximalNonbases:CongruenceNonbasis-N0-aa}
\gcd(h,s-t) = 1.
\end{equation}
Let $Y_0$ be an infinite set of  nonnegative integers with infinite gaps, and let
\[
X_0 = \N_0 \setminus Y_0. 
\]
The set  
\[ 
A_{X_0} = \{s\} \cup \{hx+t: x \in X_0 \}
\]
is an asymptotic nonbasis of order $h$ for $\N_0$ 
that is not a subset of a maximal asymptotic nonbasis of order $h$ for ${\N_0}$. 
\et

\begin{proof}
As in the proof of Theorem~\ref{MaximalNonbases:theorem:CongruenceNonbasis}, 
we begin by showing that there is a finite set $\mathcal{F}$ of nonnegative integers such that 
\begin{equation}                                                   \label{MaximalNonbases:CongruenceNonbasis-hA-N0}
\N_0 \setminus hA_{X_0} = \mathcal{F} \cup \{ (h-1)s + hy+t : y \in Y_0 \}
\end{equation} 
and 
\begin{equation}                                                   \label{MaximalNonbases:CongruenceNonbasis-hA-N0-F}
\mathcal{F} \cap \{ (h-1)s + hy+t : y \in Y_0 \} = \emptyset.
\end{equation}

Let $n$ be a nonnegative integer such that 
\begin{equation}                               \label{MaximalNonbases:CongruenceNonbasis-N0-bb}
n \not\equiv t-s \pmod{h}
\end{equation}
and
\begin{equation}                                                   \label{MaximalNonbases:CongruenceNonbasis-ineq}
n \geq (h-2) s + ht. 
\end{equation}
The divisibility condition~\eqref{MaximalNonbases:CongruenceNonbasis-N0-aa} 
and the congruence condition~\eqref{MaximalNonbases:CongruenceNonbasis-N0-bb} 
imply that there is a unique integer $i \in \{0,1,\ldots, h-2\}$ such that 
\[
n \equiv i(s-t)\pmod{h}.
\]
Let $k = h-i \in \{2,3 \ldots, h\}$.  There is a unique integer $q$ such that 
\[
n = i(s-t) + hq = is + h(q - t) +kt.
\]
Inequality~\eqref{MaximalNonbases:CongruenceNonbasis-ineq} implies that 
\[
q - t = \frac{n -is- kt}{h} \geq \frac{n - (h-2) s - ht}{h} \geq 0.
\]
By Lemma~\ref{MaximalNonbases:lemma:Ygap}, 
for all but finitely many nonnegative integers $q$, there exist $x_1,\ldots, x_k \in X_0$ such that 
\[
x_1+\cdots + x_k = q-t
\]
and so  
\begin{align*}
n & = is + h(q - t) +kt \\
& = is + h \left(x_1+\cdots + x_k \right)  +kt \\
& = is + \sum_{j=1}^k (hx_j +t) \in hA_{X_0}. 
\end{align*}   
Thus, the set $\mathcal{F}_0$ of nonnnegative integers $n$ 
such that $n \not\equiv t-s \pmod{h}$ and $n \notin hA_{\N_0}$ 
is finite.

Let $n$ be a nonnegative integer such that 
\[
n \equiv t-s \pmod{h}.
\]
By Theorem~\ref{MaximalNonbases:theorem:CongruenceBasis},  
there is a unique integer $z_1 \in \Z$ such that 
\begin{equation}                               \label{MaximalNonbases:n-rep}
n = (h-1)s + hz_1+t 
\end{equation}   
and this is the unique representation of $n$ as a sum of $h$ elements 
in the set  $A_{\Z} = \{s\} \cup \{hz+t: z\in \Z\}$.
Moreover,
\[
z_1 \geq 0 \qquad \text{if and only if} \qquad  n \geq (h-1)s +t .
\]
Let $\mathcal{F}_1$ be the finite set of nonnegative integers $n$ such that 
$n < (h-1)s +t$ and $n \equiv t-s\pmod{h}$.
Let 
\[
\mathcal{F} = \mathcal{F}_0 \cup \mathcal{F}_1.
\]
In the representation~\eqref{MaximalNonbases:n-rep}, 
if $z_1 = x_1\in X_0$, then $hx_1+t \in A_{X_0}$ and $n \in hA_{X_0}$. 
If $z_1 = y_1 \in Y_0$, then $hy_1+t \notin A_{X_0}$ and $n \notin h A_{X_0}$.  
This proves~\eqref{MaximalNonbases:CongruenceNonbasis-hA-N0}.  
Because the set $Y_0$ is infinite, 
the set $A_{X_0}$ is an asymptotic nonbasis of order $h$ for $\N_0$.

We shall prove that $A_{X_0}$ is not a subset of a maximal asymptotic nonbasis 
of order $h$ for ${\N_0}$.   

Let  $b \in  \N_0\setminus  A_{X_0}$ satisfy   
\[
b \not\equiv s,t \pmod{h}.
\]
As in the proof of Theorem~\ref{MaximalNonbases:theorem:CongruenceNonbasis}, the congruence 
\begin{equation}                                                \label{MaximalNonbases:CongruenceNonbasis-congruence-st}
(i+1)(s-t) \equiv t-b \pmod{h}.
\end{equation}
has a unique solution  $i \in \{0,1,2,\ldots, h-3 \}$.  
Let $n \in \N_0 \setminus (hA_{X_0} \cup \mathcal{F})$ satisfy 
\begin{equation}                                                \label{MaximalNonbases:CongruenceNonbasis-ineq-st}
n \geq b + (h-3)s + (h-1)t.
\end{equation}
By~\eqref{MaximalNonbases:CongruenceNonbasis-hA-N0}, there exists $y \in Y_0$ such that 
\[
n = (h-1)s+hy+t.
\]
The congruence~\eqref{MaximalNonbases:CongruenceNonbasis-congruence-st} 
and the inequality~\eqref{MaximalNonbases:CongruenceNonbasis-ineq-st} 
imply that there exists a unique nonnegative integer  $w$ such that 
\[ 
n = (h-1)s + hy + t =  b + is + hw  + (h-i-1)t. 
\]
Because $h-i-1\geq 2$,  Lemma~\ref{MaximalNonbases:lemma:Ygap} implies that, 
for all but finitely many $w$, 
there exist integers $x_1,x_2,\ldots, x_{h-i-1} \in X_0$ 
such that 
\[
w = \sum_{j=1}^{h-i-1} x_j
\]
and so 
\begin{align*}
n & =  b+ is  +  hw + (h-i-1)t \\
& =  b+ is  + h \sum_{j=1}^{h-i-1} x_j  + (h-i-1)t \\
& =  b+ is  +  \sum_{j=1}^{h-i-1} ( hx_j +t) \\
& \in  hA_{X_0}.
\end{align*} 
Thus, if $b \in  \N_0\setminus A_{X_0}$ and $b \not\equiv s,t \pmod{h}$, 
then $A_{X_0} \cup \{b\}$ is an asymptotic basis of order $h$ for $\N_0$.

Let $b \in  \N_0\setminus A_{X_0}$ satisfy 
\[
b \equiv s \pmod{h}.
\]
We have  $b \neq s$ and $b = s+hu$ for some $u \neq 0$.
For all $y \in Y_0$, we have 
\begin{align*}
(h-1)s+ hy + t  & = (h-1)(s+hu) + h(y-(h-1)u) + t \\ 
 & = (h-1)b + h(y-(h-1)u) + t.
\end{align*}
Because $Y_0$ has infinite gaps, there are only finitely many integers $y \in Y_0$ 
with $y-(h-1)u \in Y_0$, 
and so $y-(h-1)u \in X_0$ for all but finitely many $y \in Y_0$. 
Therefore, $(h-1)s+ hy + t \in h(A_{X_0} \cup \{b\})$ for all but finitely many $y \in Y_0$. 
Thus, if $b \in  \N_0\setminus A_{X_0}$  and $b \equiv s \pmod{h}$, then $A_{X_0} \cup \{b\}$ 
is an asymptotic basis of order $h$ for $\N_0$.  

We have proved that if $b \in  \N_0\setminus A_{X_0}$ and  $A_{X_0} \cup \{ b\}$ is an asymptotic nonbasis 
of order $h$ for $\N_0$, then    
\[
b \equiv t \pmod{h}
\]
and   
\[
b = hy' + t \qquad \text{for some $y' \in Y_0$.} 
\]
It follows from Theorem~\ref{MaximalNonbases:theorem:CongruenceBasis} that 
\[
h \left( A_{X_0} \cup \{b\} \right) =  hA_{X_0} \cup \mathcal{F'} \cup \{(h-1)s+b\}  
\]
for some $\mathcal{F'} \subseteq \mathcal{F}$. 
Therefore, if $B \subseteq \N_0\setminus A_{X_0}$ and $A_{X_0} \cup B$ is an asymptotic nonbasis 
of order $h$, then there exists $Y_0' \subseteq Y_0$ such that 
\[
B = B_{Y'_0} = \{hy'+t:y'  \in Y_0' \} 
\]
and 
\[
h(A_{X_0}\cup B_{Y'_0} ) = hA_{X_0}  \cup \mathcal{F''}  \cup \{ (h-1)s+hy' +t:y' \in  Y_0'  \}.  
\] 
for some $\mathcal{F''} \subseteq \mathcal{F}$. 
We conclude that $A_{X_0}\cup B_{Y'_0} $ is an asymptotic nonbasis of order $h$ 
for ${\N_0}$\ if and only if $Y_0'$ is a co-infinite subset of $Y_0$.
Because the infinite set $Y_0$ contains no maximal co-infinite subset, 
the set $A_{X_0}$ is an asymptotic nonbasis of order $h$ 
for ${\N_0}$\ that is not contained
in a maximal asymptotic nonbasis of order $h$ for ${\N_0}$.
This completes the proof.  
\end{proof}

\def\cprime{$'$} \def\cprime{$'$} \def\cprime{$'$}
\providecommand{\bysame}{\leavevmode\hbox to3em{\hrulefill}\thinspace}
\providecommand{\MR}{\relax\ifhmode\unskip\space\fi MR }
\providecommand{\MRhref}[2]{%
  \href{http://www.ams.org/mathscinet-getitem?mr=#1}{#2}
}
\providecommand{\href}[2]{#2}

\end{document}